\documentclass{article}
\usepackage{amssymb}
\usepackage{amsmath}
\usepackage{amsthm}
\newcommand{\qq}{\mathbb Q} \newcommand{\nn}{\mathbb N} 
\newcommand{\zz}{\mathbb Z} \newcommand{\rr}{\mathbb R}   
\newtheorem{thm}{Theorem}

\newtheorem{exam}{Example}

\begin{document}
\Large

\title{\LARGE\bf Note on Archimedean property in ordered vector spaces}
\author{\bf Eduard Yu. Emelyanov\footnote{Middle East Technical University, Ankara, Turkey}}

\date{}
\maketitle

{\bf Abstract:} {\normalsize 
It is shown that an ordered vector space $X$ is Archimedean if and only if
$\inf\limits_{\tau\in\{\tau\}, y\in L}(x_\tau -y) \ = 0$ for any bounded decreasing net $x_\tau\downarrow$ in $X$, 
where $L$ is the collection of all lower bounds of $\{x_\tau\}_{\tau}$. 
We give also a characterization of the almost Archimedean property of $X$ in terms of existence of a linear  
extension of an additive mapping $T:Y_+\to X_+$ of the positive cone $Y_+$ of an ordered vector space $Y$
into $X_+$.}



\section{Introduction}

In this note we deal with vector spaces over reals. 
A subset $K$ of a vector space $X$ is called a {\bf cone} if it satisfies 
$$
  K\cap(-K)=\{0\} \text{,} \ \ K+K\subseteq K  \ \ \text{\rm and} \ \ rK\subseteq K
$$ 
for all $r\ge 0$. A cone $K$ is said to be {\bf generating} if $K-K=X$.
Given a cone $X_+$ in $X$, we say that $(X,X_+)$ is an {\bf ordered vector space}.
The {\bf partial ordering} $\le$ on $X$ is defined by 
$$ 
  x\le y \ \ \text{\rm if} \ \ y-x\in X_+\,.
$$ 
The space $(X,X_+)$ is also denoted by $(X,\le)$ or simply by $X$ if $\le$ is well understood.

In what follows, we denote by $X$ an ordered vector space.
For every $x,y\in X$, the (possibly empty if $x\not\le y$) set 
$$
  [x,y]=\{z:x\le z\le y\}
$$ 
is called {\bf order interval}. 
A subset $A\subseteq X$ is said to be {\bf order convex} 
if for all $a,b\in A$, we have $[a,b]\subseteq A$.
Given $B\subseteq X$, the smallest order convex subset $[B]$ of $X$ containing $B$ is called the 
{\bf order convex hull} of B. It is immediate to see that 
$$
  [B]=(B+X_+)\cap(B-X_+)=\bigcup_{a,b\in B}[a,b].
$$
Any order convex vector subspace of $X$ is called an {\bf order ideal}. In the case when 
$Y\subseteq X$ is an order ideal, the quotient space $X/Y$ is partially ordered by:
$$
  [0]\le [f] \ \ \text{\rm if} \ \ \exists q\in Y \ \text{\rm with} \ \  0\le f+q\ .
$$
The order convexity of $Y$ is needed for the property:
$$
  [f]\le [g] \le [f] \ \Rightarrow [f]=[g]\,,
$$
which guarantees that the canonical image $X_+$ is a cone in  $X/Y$.
Indeed if $[f]\le [g] \le [f]$, then 
$$
  q_1\le g-f \ \ \text{\rm and} \ \ -q_2\le f-g
$$
for some $q_1,q_2\in Y$. Thus, $q_1\le g-f\le q_2$. Since $Y$ is order convex, then $g-f\in [q_1,q_2]\subseteq Y$ and $[g]=[f]$.

Let $(x_n)$ be a sequence in $X$ and $u\in X_+$. The sequence $(x_n)$ is said to be {\bf $u$-uniformly convergent}
to a vector $x\in X$, in symbols $x_n\stackrel{(u)}{\to}x$, if 
$$
  x_n-x\in[-\varepsilon_nu,\varepsilon_nu]
$$ 
for some sequence $(\varepsilon_n)$ of reals such that $\varepsilon_n\downarrow 0$. A subset $A\subseteq X$ is said to be 
{\bf uniformly closed} in $X$ if it contains all limits of all $u$-uniformly convergent sequences $(x_n)\subseteq A$ for all $u\in X_+$.

An ordered vector space $X$ is said to be {\bf Archi\-me\-dean} (we say also that $X_+$ has the {\bf Archi\-medean property}) if
$$
  [(\forall n\ge 1) \ ny\le x\in X_+]\Rightarrow [y\le 0]\,.
  \eqno(1)
$$
It is easy to see that in $(1)$ a vector $x\in X_+$ can be replaced by any $x\in X$. 
Any subspace of an Archimedean ordered vector space is Archimedean. It is well known 
that $X$ is Archimedean iff $\inf_{n\ge 1}\frac{1}{n}y=0$ for every $y\in X_+$.
It is worth to remark that if $X$ admits a linear topology $\tau$ for which its cone $X_+$ is closed, 
then $X$ is Archimedean. Indeed, assume $ny\le x$ for all $n\ge 1$ and some $x,y\in X$. Then 
$$
  \frac{1}{n}x-y\in X_+ \ \ \text{\rm and} \ \ \frac{1}{n}x-y\stackrel{\tau}{\to} -y
$$ 
imply $-y\in X_+$ or $y\le 0$.

$X$ is called {\bf almost Archimedean} if
$$
  [(\forall n\in\zz) \ ny\le x\in X_+]\Rightarrow [y=0]\,,
$$
for every $x\in X_+$. Clearly, $X$ is almost Archimedean iff
$$
  \bigcap_{n\ge 1}\Bigl[-\frac{1}{n}x,\frac{1}{n}x\Bigr]=\{0\}\,,
$$
for every $x\in X_+$. It follows immediately, that any subcone of an almost Archimedean cone is almost Archi\-medean.  
A standard example of an ordered vector space which is not almost Archimedean is $(\rr^2,\le_{lex})$, 
the Euclidean plain with the lexicographic ordering. If $Y\subseteq X$ is an order ideal such that $X/Y$ is almost Archimedean, 
then $Y$ is uniformly closed. Indeed, let $y_n\stackrel{(u)}{\to}x$ for some $(y_n)\subseteq Y$, $u\in X_+$, and $x\in X$. 
We may assume that
$$
  y_n-x\in\Bigl[-\frac{1}{n}u,\frac{1}{n}u\Bigr] \ \ \ (\forall n\ge 1).
$$ 
Hence
$$
  [-x]=[y_n-x]\in \Bigl[-\frac{1}{n}[u],\frac{1}{n}[u]\Bigr] \ \ \ (\forall n\ge 1).
$$ 
Since $X/Y$ is almost Archimedean, then $[-x]=[0]$, and thus $x\in Y$. If $Y=\{0\}$, the converse is obviously true. 
Thus $X$ is almost Archimedean iff $\{0\}$ is uniformly closed. In general, the question whether or not $X/Y$ is almost Archimedean
assuming an order ideal $Y$ to be uniformly closed is rather nontrivial (it has a positive answer, for example in vector lattice setting 
(cf. \cite[Thm.60.2]{LZ})). Any Archimedean ordered vector space is clearly almost Archi\-medean. The converse is not true even when $\dim(X)=2$.
\begin{exam}
{\it 
Let $\Gamma$ a set containing at least two elements and let $Y$ be the space of all bounded real functions
on $\Gamma$, partially ordered by:
$$
  f\le g \ \ \text{\it if} \ \ \text{\it either} \ f=g \ \  \text{\it or} \ \ \inf_{t\in\Gamma}[g(t)-f(t)]>0\,.
$$
The space $(Y,\le)$ is almost Archi\-me\-dean but not Archi\-me\-dean. Indeed, it can be seen easily that
$\inf_{n\ge 1}\frac{1}{n}f$ does not exist for any $0\ne f\in Y_+$.  
}
\end{exam}

An ordered vector space $X$ is said to be a {\bf vector lattice} (or a {\bf Riesz space}) if every nonempty finite subset of $X$ has a 
least upper bound. $X_+$ is said to be {\bf minihedral} if $X$ is a vector lattice. Any almost Archimedean vector lattice is Archimedean,
indeed
$$
  \Bigl[(\forall n\ge 1)(ny\le x\in X_+)\Bigr]\Rightarrow 
$$
$$
  \Bigl[(\forall n\ge 1)[-x\le n\sup(y,0)\le \sup(x,0)=x]\Bigr]\Rightarrow
$$
$$
  \Bigl[\sup(y,0)=0\Bigr]\Rightarrow \Bigl[y\le 0\Bigr]\,.
$$

For further details on ordered vector spaces we refer to the book \cite{AT}.

\section{A characterization of Archimedean ordered vector spaces}

The following characterization of the Archimedean property is well known in the vector lattice case 
(see, for example, \cite[Thm.22.5]{LZ}). In the general setting of ordered vector spaces, it has
been published recently \cite[Prop.2]{E} as an auxiliary fact with an incorrect proof of 
the implication $(b)\Rightarrow(a)$. Below, we fill the gap in the proof. 

\begin{thm}
{\it For an ordered vector space $X$, the following are equivalent$:$

$(a)$  \  $X$ is Archimedean.

$(b)$  \  For any decreasing net $x_\tau\downarrow\ge d$ in $X$,
$$
  \inf\limits_{\tau\in\{\tau\}, y\in L}(x_\tau -y) \ = 0\,,
$$
where $L=\{y\in X: (\forall \tau\in\{\tau\})[y\le x_\tau]\}$ is the collection of all lower bounds of
$\{x_\tau\}_{\tau}$.}
\end{thm}

{\bf Proof}:  $(a)\Rightarrow(b)$:  Let $X$ be Archimedean, $x_\tau\downarrow \ge d$. 
Assume $z\in X$ satisfies $z\le x_\tau - y$ for all indexes $\tau$ and for all $y\in L$. 
Since $0\le x_\tau - y$ for all $\tau$, to complete the proof of the implication, it is enough to show that $z\le 0$.

As $y+z\le x_\tau$ for all $\tau$ and all $y\in L$, we obtain that $y+z\in L$ for every  $y\in L$. 
It follows by the induction, 
$$
  y+nz\in L \ \ \ (\forall y\in L)(\forall n\in\nn)\,.
$$ 
In particular, $d+nz\le x_{\tau_0}$ (and hence, $nz\le x_{\tau_0}- d$) 
for some $\tau_0\in\{\tau\}$ and all $n\in\mathbb{N}$. By the condition 
$$
  nz\le x_{\tau_0}- d \ \ (\forall n\in\mathbb{N})\,,
$$
the Archimedean property of $X_+$ implies $z\le 0$, what is required.

$(b)\Rightarrow(a)$:  Let $x\in X_+$, $ny\le x$ for all $n\in\mathbb{N}$. We have to show that $y\le 0$.
Denote 
$$
  L=\Bigl\{w\in X: (\forall n\ge 1)\ w\le \frac{1}{n}x\Bigr\}\,.
$$
Clearly, $y\in L$. Given $u\in L$, then $0\le \frac{1}{n}x-u$ for all $n\ge 1$. By the hypothesis applied to the decreasing sequence 
$\frac{1}{n}x\downarrow\ge y$, the following infima exist, and
$$
  \inf\limits_{n\ge 1, u\in L}\Bigl(\frac{2}{n}x-u\Bigr)=\inf\limits_{n\ge 1, u\in L}\Bigl(\frac{1}{n}x-u\Bigr)=0\,.
$$
Hence,
$$  
  \inf\limits_{n\ge 1, u\in L}\Bigl[\Bigl(\frac{2}{n}x-u\Bigr)-y\Bigr]=-y+\inf\limits_{n\ge 1, u\in L}\Bigl(\frac{2}{n}x-u\Bigr)=-y\,. 
$$
Since
$$
  0\le \Bigl(\frac{1}{n}x-u\Bigr)+\Bigl(\frac{1}{n}x-y\Bigr)=\Bigl(\frac{2}{n}x-u\Bigr)-y 
$$
for all $n\ge 1$, then 
$$
  0\le \inf\limits_{n\ge 1, u\in L}\Bigl[\Bigl(\frac{2}{n}x-u\Bigr)-y\Bigr]=-y\,.
$$
Then, $y\le 0$ and hence $X$ is Archimedean. $\blacksquare$

\section{A characterization of almost Archimedean ordered vector spaces}

Here we present a characterization of the almost Archi\-medean property of $X_+$ in terms of existence of an 
extension of an additive mapping $T:Y_+\to X_+$ of the positive cone $Y_+$ of an ordered vector space $Y$
into $X_+$ to a linear operator from $Y$ into $X$. The existence of such an extension is well known 
in the Archimedean setting (see, for example, \cite[Lem.1.26]{AT}, \cite[Lem.83.1]{Z}). 

\begin{thm}
{\it 
For an ordered vector space $X$, the following statements are equivalent$:$

$(i)$  \ \ $X$ is almost Archimedean.

$(ii)$  \  For any ordered vector space $Y$ and any additive mapping $T:Y_+\to X_+$,
there is an extension of $T$ to a linear operator from $Y$ to $X$.

$(iii)$    For any additive mapping $T:\rr_+\to X_+$, there is an extension of $T$
to a linear mapping from $\rr$ to $X$.
}
\end{thm}

{\bf Proof}: $(i)\Rightarrow(ii)$: Since on an algebraic complement $Y_0$ of $Y_+-Y_+$ in $Y$, an extension of $T$, say $S$, can be chosen as any linear operator, 
we only need to obtain an extension of an additive mapping $Y_+\stackrel{T}{\to}X_+$ to a linear operator $(Y_+-Y_+)\stackrel{S}{\to}X$.

For each $y\in (Y_+-Y_+)$ pick $y_1,y_2\in Y_+$ with $y=y_1-y_2$ and put 
$$
  Sy=T(y_1)-T(y_2)\,.
$$
It is routine to see that $S$ is well defined and additive. 
The additivity of $(Y_+-Y_+)\stackrel{S}{\to}X$ implies that $S$ is $\qq$-homogeneous. 
To complete the proof, it is enough to show that $S$ is $\rr_+$-homogeneous on $Y_+$, where it coincides with $T$.
Thus we need to show that $T:Y_+\to X_+$ is $\rr_+$-homogeneous. We shall use the following elementary remark:
$$
  \Bigl[x,y\in [q,p] \Bigr] \Rightarrow \Bigl[x-q,y-q\in [0,p-q]\Bigr] \Rightarrow
$$
$$  
  \Bigl[x-y=(x-q)-(y-q)\in [-(p-q),p-q]\Bigr].
  \eqno(2)
$$
Let $\qq_+\ni r_n\uparrow a\in\rr$, $\qq_+\ni r_n'\downarrow a$, $u\in Y_+$. Then
$$
  r_nT(u)=T(r_nu)\le T(au)\le T(r_n'u)=r_n'T(u)\,.  
$$
But also:
$$
  r_nT(u)\le aT(u)\le r_n'T(u)\,.  
$$
That is 
$$
  T(au),aT(u)\in [r_nT(u),r_n'T(u)]\,.
  \eqno(3)
$$
Applying $(2)$ in $(3)$, we obtain 
$$
  T(au)-aT(u)\in [-(r_n'-r_n)Tu,(r_n'-r_n)Tu] \ \ (\forall n\in\nn)\,.
$$
Since $X$ is almost Archimedean and $r_n'-r_n\downarrow 0$, 
we get $T(au)-aT(u)=0$, what is required.

$(ii)\Rightarrow(iii)$: It is trivial.

$(iii)\Rightarrow(i)$: Let $x\in X_+$, $y\in X$ be such that
$$
  -\frac{1}{n}x\le y\le \frac{1}{n}x \ \ \ (\forall n\ge 1).
  \eqno(4)
$$
Take a function $f:\rr\to\rr$ which is $\qq$-linear but not $\rr$-linear and define an additive mapping $T:\rr_+\to X$ by
$$
  T(a)=ax+f(a)y \ \ \ (a\in \rr_+)\,.
$$
Then $T$ maps $\rr_+$ into $X_+$. Indeed, $T(0)=0$ and if $0<a$ then $0\le a-\frac{|f(a)|}{n_0}$ for some large enough $n_0$.
It follows from (4) that $-\frac{|f(a)|}{n}x\le \pm f(a)y\le \frac{|f(a)|}{n}x$ for all $n\ge 1$.
In particular, $-\frac{|f(a)|}{n_0}x\le f(a)y$, and therefore
$$
  0\le \Bigl(a-\frac{|f(a)|}{n_0}\Bigr)x\le ax+f(a)y = T(a)
$$
Take a linear extension $S$ of $T$ to all of $\rr$. Then $S$ (and hence $T$) must be $\rr_+$-homogeneous on $\rr_+$ which is only possible if $y=0$. 
$\blacksquare$

\end{document}